\documentclass[12pt,reqno]{amsart}
\usepackage[usenames]{color}
\usepackage{amsmath, amsfonts, epsfig, amssymb, graphics}
\usepackage[colorlinks=true,citecolor=cyan]{hyperref}
\usepackage[mathscr]{euscript}

\newcommand{\pref}[1]{{\rm (\ref{#1})}}
\def\mb{\mathbf}

\def\dinv{\mathrm{dinv}}
\def\co{\bleu{\mathbf{co}}}

\def\R{\mathcal{R}}
\def\Altid{\mathcal{A}}
\def\Alt{\mathscr{A}}
\def\Ideal{\mathcal{I}}
\def\Inv{\mathscr{I}}
\def\Base{\mathscr{B}}
\def\C{\mathbb{C}}
\def\Harm{\mathscr{H}}

\def\S{{\mathbb S}}
\def\N{{\mathbb N}}

\newcommand{\sign}{{\rm sign}}

\newcommand{\scalar}[2]{{\langle#1,#2 \rangle}}

\def\auteur#1{{\sc #1}}
\def\titreref#1{{\em #1}}
\def\vol#1{{\bf #1}}
\def\defn#1{\bleu{{\bf #1}}}

\newcommand{\del}[1]{\partial#1}
\newcommand{\Park}{\mathcal{P\!F}}
\newcommand{\Parking}[1]{\Park^{(r)}\!(#1)}
\def\Dyck#1#2{\mathcal{D}^{(#1)}_{#2}}
\def\charac{\raise 2pt\hbox{\begin{math}\chi\end{math}}}

\newdimen\carrelength

\def\bleu{\textcolor{blue}}
\def\rouge{\textcolor{red}}
\def\jaune{\textcolor{yellow}}

\newtheorem*{thm}{\bleu{Theorem}}
\newtheorem{proposition}{\bleu{Proposition}}
\newtheorem{lemma}{\bleu{Lemma}}
\newtheorem{conjecture}{\bleu{Conjecture}}

\begin{document}

\title[Trivariate Diagonal Harmonics]{\bleu{Higher Trivariate Diagonal Harmonics via generalized Tamari Posets}}
\author[F.~Bergeron and L.-F. Pr\'eville-Ratelle]{F. Bergeron and L.-F. Pr\'eville-Ratelle}
\address{D\'epartement de Math\'ematiques, UQAM,  C.P. 8888, Succ. Centre-Ville, 
 Montr\'eal,  H3C 3P8, Canada.}\date{May 2010}
\maketitle
\begin{abstract}
We consider the graded $\S_n$-modules of higher diagonally harmonic polynomials in three sets of variables (the trivariate case), and show that they have interesting ties with generalizations of the Tamari poset and parking functions. In particular we get several nice formulas for the associated Hilbert series and graded Frobenius characteristics. This also leads to entirely new combinatorial formulas.
\end{abstract}
 \parskip=0pt

{ \setcounter{tocdepth}{1}\parskip=0pt\footnotesize \tableofcontents}
\parskip=8pt

\section{Introduction}
Our global aim here is to give an explicit description of the space $\bleu{\Harm_n}$ of \defn{trivariate diagonal harmonic polynomials}. We make interesting progress along these lines, and propose a combinatorial description of the graded Frobenius characteristic for this space. This description encompasses the ``shuffle conjecture'' of~\cite{HHLRU}, and suggests new combinatorial identities linked to the combinatorics of the Tamari lattice. In fact it has lead us to introduce a new notion of $r$-Tamari poset, and these identities may be stated in terms of positive integer polynomials in the variable $r$.
The first drafts of this paper has prompted other researchers to consider these combinatorial questions, and some have already been settled (see~\cite{melou} and \cite{chapuy}). 
  
More explicitly, in the  ring $\bleu{\R_n:=\C[X]}$ of polynomials in the three sets of variables
   \begin{displaymath}\bleu{X:=\begin{pmatrix}
      x_1 & x_2 & \cdots &x_n\\
      y_1& y_2 & \cdots &y_n\\
      z_1&  z_2 & \cdots &z_n
      \end{pmatrix}, }\end{displaymath}
  ${\Harm_n}$ is the subspace of polynomial zeros of the 
\defn{polarized power sums} differential operators
  \begin{equation}\label{def_power_sum}
      \bleu{P_\alpha(\partial X):=\sum_{j=1}^n\del{x_j}^a\del{y_j}^b \del{z_j}^c},
  \end{equation}
for $\alpha=(a,b,c)$ running thought the set of vectors of \defn{norm} $\bleu{|\alpha|:=a+b+c}$, with $1\leq |\alpha|\leq n$. Here,  $\partial v$ denotes derivation with respect to $v$. It is clear from the definition that $\Harm_n$ is closed under derivation. 
For one set of variables, say $\mb{x}=x_1,\ldots,x_n$,  $\Harm_n$ is entirely described by classical theorems (see~\cite[section 3.6]{humphreys}) as the span of all partial derivatives of the Vandermonde determinant
     \begin{displaymath}\Delta_n(\mb{x})=\det\left( x_i^j\right)_{1\leq i\leq n,\ 0\leq j\leq n-1},\end{displaymath}
It is of dimension $n!$, and isomorphic to the regular representation of $\S_n$, for the action that permutes variables. The description of the case of two sets of variables took close to 15 years to be finalized (see~\cite{haimanhilb}), and much remains to be understood. It is of dimension $(n+1)^{n-1}$. It has been observed by Haiman (see \cite{JAC}) that the dimension of $\Harm_n$ seems to be given by the formula
\begin{equation}\label{dim_Hn}
   \bleu{\dim \Harm_n= 2^n\,(n+1)^{n-2}},
\end{equation}
but, until very recently, almost no one had further studied the trivariate case. We have endeavoured to do so, considering more general spaces $\bleu{\Harm_{n}^{(r)}}$ (see definition~\pref{def_higher}) for which we have experimentally found out that 
\begin{equation}\label{dim_Hrn}
   \bleu{\dim\Harm_{n}^{(r)}= (r+1)^n\,(r\,n+1)^{n-2}}.
\end{equation}

\subsection*{Symmetric group action}
We turn both $\R_n$ and $\Harm_n$ into $\S_n$-module  by considering the \defn{diagonal action} of symmetric group $\S_n$ on  variables.  Recall that, for $\sigma\in \S_N$, the polynomial $\sigma\cdot f(X)$ is obtained by replacing the variables in $f(X)$ by
     \begin{displaymath}\begin{pmatrix}
      x_{\sigma(1)} & x_{\sigma(2)} & \cdots &x_{\sigma(n)}\\
      y_{\sigma(1)}& y_{\sigma(2)} & \cdots &y_{\sigma(n)}\\
      z_{\sigma(1)}&  z_{\sigma(2)} & \cdots &z_{\sigma(n)}
      \end{pmatrix}. \end{displaymath}
We denote by  $X^A$ the monomial
    \begin{displaymath}\bleu{X^A:=x_1^{a_1}\cdots x_n^{a_n}\cdot y_1^{b_1}\cdots y_n^{b_n}\cdot z_1^{c_1}\cdots z_n^{c_n}},\end{displaymath} 
 for a matrix of integers
    \begin{displaymath}A=\begin{pmatrix}
      a_1 & a_2 & \cdots &a_n\\
      b_1& b_2 & \cdots &b_n\\
      c_1&  c_2 & \cdots &c_n
      \end{pmatrix}. \end{displaymath}
Writing $A_j$ for the  $j^{\rm th}$-column of $A$, we may define
the \defn{degree} of a monomial $X^A$ to be the vector 
    \begin{displaymath}\bleu{\deg(X^A):=\sum_{j=1}^n A_j},\end{displaymath} 
 in $\N^3$. Clearly the diagonal action preserves degree.
 The \defn{total degree} ${\rm tdeg}(X^A)$ of a monomial $X^A$ is the sum of the components of $\deg(X^A)$.
An interesting subspace of $\Harm_n$ is the space $\bleu{\Alt_n}$ af \defn{alternating} polynomials in $\Harm_n$. Recall that these are the polynomials such that $\sigma\cdot f(X)=\sign(\sigma)\,f(X)$. It is also worth recalling that, in the case of two sets of variables, this subspace has dimension equal to the Catalan number
   \begin{displaymath}
         \bleu{C_n:=\frac{1}{n+1}\binom{2\,n}{n}}.
    \end{displaymath}
 Another observation of Haiman ({\sl loc. sit.}) is that, in the trivariate case, we seem to have
\begin{equation}\label{dim_An}
   \bleu{\dim\Alt_n=\frac {2}{ n\,( n+1) } \binom{4\,n+1}{ n-1}}.
\end{equation}
More generally, our experimental calculations suggest that
 the  alternating component $ \Alt_n^{(r)}$ of $ \Harm_n^{(r)}$ has dimension
\begin{equation}\label{dim_Altr}
   \bleu{ \dim \Alt_n^{(r)}}=\bleu{\frac {( r+1)}{ n\,( r n+1) } \binom{( r+1) ^{2}\,n+r}{ n-1} }.
\end{equation}
It is intriguing\footnote{Thanks to  the anonymous referee for this observation.} that this formula bears a resemblance to a formula that appears in \cite{kontsevich}.

\subsection*{Degree grading} The space $\Harm_n$ is graduated by degree, and its \defn{homogeneous components}  $\bleu{\Harm_{n,\mb{d}}}$ (with $d\in \N^3$) are $S_n$-invariant, hence ${\Alt_n}$ is also graded. There corresponds direct sum decompositions
\begin{equation}\label{decomp_hom}
   \bleu{ \Harm_n=\bigoplus_{d\in \N^3} \Harm_{n,\mb{d}}},\qquad {\rm and}\qquad  \bleu{ \Alt_n=\bigoplus_{d\in \N^3} \Alt_{n,\mb{d}}}.
 \end{equation} 
As discussed in Section~\ref{invariant}, the space $\Harm_n$ is finite dimensional. In fact, we will see that the homogeneous component $\Harm_{n,\mb{d}}$ is non vanishing only if 
    \begin{displaymath}|\mb{d}|=d_1+d_2+d_3\leq \binom{n}{2},\end{displaymath}
where $\mb{d}=(d_1,d_2,d_3)$.  We can thus consider the \defn{Hilbert series}:
   \begin{equation}\label{defn_Hilb}
        \bleu{\Harm_n(\mb{q}):=\sum_{|\mb{d}|=0}^{\binom{n}{2}} \dim(\Harm_{n,\mb{d}})\, \mb{q}^\mb{d}}, \qquad{\rm and}\qquad
         \bleu{\Alt_n(\mb{q}):=\sum_{|\mb{d}|=0}^{\binom{n}{2}} \dim(\Alt_{n,\mb{d}})\, \mb{q}^\mb{d}}
    \end{equation}
  of $\Harm_n$ and $\Alt_n$,  writing $\mb{q}^\mb{d}$ for $q_1^{d_1}q_2^{d_2}q_3^{d_3}$. The Hilbert series $\Harm_n(\mb{q})$ is a symmetric polynomial in $q_1$, $q_2$ and $q_3$. In fact, it is always Schur positive as discussed in Section~\ref{general}. This is to say that it expands as a positive integer coefficient linear combination of the Schur polynomials $s_\mu(\mb{q})$ in the variables $\mb{q}$. 
 
For example, we easily calculate that $\Harm_2$ affords the linear basis
\begin{equation}\label{base_2}
   \{1,x_2-x_1,y_2-y_1,z_2-z_1\},
 \end{equation}
whose subset $ \{x_2-x_1,y_2-y_1,z_2-z_1\}$ is clearly a basis of  $\Alt_2$.
We thus get  the Hilbert series
   \begin{displaymath}\Harm_2(\mb{q})=1+(q_1+q_2+q_3)=1+s_{1}(\mb{q}),\qquad {\rm and}\qquad \Alt_2(\mb{q})=q_1+q_2+q_3=s_{1}(\mb{q}).\end{displaymath}
To calculate $\Harm_n(\mb{q})$ for larger $n$, we exploit the fact $\Harm_n$ is closed under taking derivatives, as well as applying the operators (see Section~\ref{invariant})
   \begin{displaymath}\bleu{E_{\mb{u}\mb{v}}^{(k)}:= \sum_{i=1}^n u_i\del{v_i}^k},\qquad (E_{\mb{u}\mb{v}}:=E_{\mb{u}\mb{v}}^{(1)}). \end{displaymath}
Here $\mb{u}$, $\mb{v}$ stand for any two of the three sets of $n$ variables:
    \begin{displaymath}\mb{x}=x_1,\ldots,x_n,\qquad \mb{y}=y_1,\ldots,y_n,\qquad {\rm and}\qquad \mb{z}=z_1,\ldots,z_n.\end{displaymath}
It easy to verify (see Section~\ref{invariant}) that the Vandermonde determinant $\Delta_n(\mb{u})$ 
belongs to $\Harm_n$, whether $\mb{u}$ be $\mb{x}$, $\mb{y}$ of $\mb{z}$. 
Using all this, we can calculate that
  \begin{eqnarray*}
      \Harm_3(\mb{q})&=&1+2\,(q_1+q_2+q_3)\\ 
         &&+2\,(q_1^{2}+q_2^{2}+q_3^{2})+3\,(q_1q_2+q_1q_3+q_2q_3)\\
         &&+q_1^{3}+q_2^{3}+q_3^{3}+q_1^{2}q_2+q_1^{2}q_3+q_1q_2^{2}+q_1q_3^{2}+q_2^{2}q_3+q_2q_3^{2}+
q_1q_2q_3
\end{eqnarray*}
by checking directly that we have the following respective bases $\mathcal{B}_{d}$ for each $\Harm_{n,\mb{d}}$:
\begin{displaymath}\begin{array}{lllll}
\Base_{300}=\{\Delta_3(\mb{x})\}, \qquad   \Base_{200}=\{\del_{x_1}\Delta_3(\mb{x}),\del_{x_2}\Delta_3(\mb{x})\},&
 \Base_{100}=\{\del_{x_1}^2\Delta_3(\mb{x}),\del_{x_1}\del_{x_2}\Delta_3(\mb{x})\}, \\[4pt]
\Base_{030}=\{\Delta_3(\mb{y})\}, \qquad    \Base_{020}=\{\del_{y_1}\Delta_3(\mb{y}),\del_{y_2}\Delta_3(\mb{y})\},&
 \Base_{010}=\{\del_{y_1}^2\Delta_3(\mb{y}),\del_{y_1}\del_{y_2}\Delta_3(\mb{y})\}, \\[4pt]
\Base_{003}=\{\Delta_3(\mb{z})\}, \qquad    \Base_{002}=\{\del_{z_1}\Delta_3(\mb{z}),\del_{z_2}\Delta_3(\mb{z})\},&
 \Base_{001}=\{\del_{z_1}^2\Delta_3(\mb{z}),\del_{x_1}\del_{z_2}\Delta_3(\mb{z})\}, \\[4pt]
\Base_{210}=\{E_{\mb{y}\mb{x}} \Delta_3(\mb{x})\}, \qquad  \Base_{120}=\{E_{\mb{y}\mb{x}}E_{\mb{y}\mb{x}} \Delta_3(\mb{x})\}, &
\Base_{201}=\{E_{\mb{z}\mb{x}} \Delta_3(\mb{x})\},   \\[4pt]
 \Base_{102}=\{E_{\mb{z}\mb{x}}E_{\mb{z}\mb{x}} \Delta_3(\mb{x})\}, \qquad 
 \Base_{021}=\{E_{\mb{z}\mb{y}} \Delta_3(\mb{y})\}, & \Base_{012}=\{E_{\mb{z}\mb{y}}E_{\mb{z}\mb{y}} \Delta_3(\mb{y})\},\\[4pt]
\Base_{110}=\{E_{\mb{y}\mb{x}}\del_{x_1}\Delta_3(\mb{x}),E_{\mb{y}\mb{x}}\del_{x_2}\Delta_3(\mb{x}), E_{\mb{y}\mb{x}}^{(2)} \Delta_3(\mb{x})\},\\[4pt] 
\Base_{101}=\{E_{\mb{z}\mb{x}}\del_{x_1}\Delta_3(\mb{x}),E_{\mb{z}\mb{x}}\del_{x_2}\Delta_3(\mb{x}), E_{\mb{z}\mb{x}}^{(2)} \Delta_3(\mb{x})\},\\[4pt] 
\Base_{011}=\{E_{\mb{z}\mb{y}}\del_{y_1}\Delta_3(\mb{y}),E_{\mb{z}\mb{y}}\del_{y_2}\Delta_3(\mb{y}), E_{\mb{z}\mb{y}}^{(2)} \Delta_3(\mb{y})\},\\[4pt] 
\Base_{000}=\{1\},\qquad \Base_{111}=\{E_{\mb{z}\mb{x}}E_{\mb{y}\mb{x}}\Delta_3(\mb{x})\}
\end{array}\end{displaymath}
Collecting previous observations, and doing some further explicit calculations, we get that
\begin{equation}\label{valeurs_Hilbert}
  \begin{array}{rcl}
 \Harm_{1}(\mb{q})&=&1,\\[3pt]
\Harm_{2}(\mb{q})&=&1+s_{1}(\mb{q})\\[3pt]
\Harm_{{3}}(\mb{q})&=&1+2\,s_{1}(\mb{q})+2\,s_{2}(\mb{q})+s_{{11}}(\mb{q})+
s_{{3}}(\mb{q}),\\[3pt]
\Harm_{{4}}(\mb{q})&=&1+3\,s_{1}(\mb{q})+5\,s_{2}(\mb{q})+3\,s_{{11
}}(\mb{q})+6\,s_{{3}}(\mb{q})+5\,s_{{21}}(\mb{q})+s_{{111}}(\mb{q})\\
    &&\qquad +5\,s_{{4}}(\mb{q})+4\,s_{{31}}(\mb{q})+3\,s_{{
5}}(\mb{q})+s_{{41}}(\mb{q})+s_{{6}}(\mb{q}).
\end{array}
\end{equation}
If we specialize one of the parameters (say $q_3$) to $0$, we get back the graded Hilbert series of bivariate diagonal harmonics (of overall dimension $(n+1)^{n-1}$) which has received a lot of attention in recent years (see \cite{cherednik, etingof, gordon, griffeth, HHLRU, haimanvanishing, loktev}). In other words, the $\mb{z}$-free (or $\mb{x}$-free, or $\mb{y}$-free) component of $\Harm_n$ coincides with the ``usual'' space of diagonal harmonics (in the bivariate case, as previously considered in the literature).
Hence, all of our formulas involving the parameters $\mb{q}$ specialize to known formulas by the simple device of setting one of the three parameters in $\mb{q}$ equal to $0$.
Evaluating $\Harm_n(\mb{q})$ at $q_i$ equal to $1$, we clearly get the overall dimension of $\Harm_n$. For this evaluation, the value of~\pref{valeurs_Hilbert} agrees with  formula~\pref{dim_Hn}.

\subsection*{Graded character} We refine the dimension analysis by taking into account the decomposition into irreducibles of the homogeneous components of $\Harm_n$, see \pref{decomp_hom}. This is all encompassed into the \defn{graded Frobenius characteristic} of $\Harm$:
    \begin{displaymath}\bleu{\Harm_n(\mb{q};\mb{w}):=\sum_{d\in \N^3} \mb{q}^\mb{d} \mathcal{F}_{\Harm_{n;\mb{d}}}(\mb{w})}.\end{displaymath}
Recall that the Frobenius characteristic $\mathcal{F}_\mathcal{V}(\mb{w})$, of a $\S_n$-module $\mathcal{V}$, is the symmetric function (in auxiliary variables $\mb{w}=w_1,w_2,\ldots$) whose expansion in terms of the Schur functions $S_\lambda(\mb{w})$ records the multiplicity of irreducibles in $\mathcal{V}$. This is to say that we have
     \begin{displaymath} \mathcal{F}_\mathcal{V}(\mb{w})= \sum_{\lambda\vdash n} n_\lambda\, S_\lambda(\mb{w}),\end{displaymath}
  with the sum being over partitions of $n$, which classify  irreducible representations of $\S_n$. For reasons also discussed in~\cite{bergeron_several}, the expansion of $\Harm_n(\mb{q};\mb{w})$ in terms of the $S_\lambda(\mb{w})$ affords Schur positive coefficients in the $s_\mu(\mb{q})$. This is to say that
   \begin{equation}
     \bleu{ \Harm_n(\mb{q};\mb{w})=\sum_{\lambda\vdash n} \left(\sum_{\mu} n_{\lambda,\mu} s_\mu(\mb{q})\right) S_\lambda(\mb{w})},\qquad n_{\lambda,\mu}\in\N.
   \end{equation}
Hence there are two kinds of Schur function playing a role here. To emphasize this, we denote those in the $\mb{q}$-variables by a lower case ``s'', and drop the variables $\mb{q}$. For example, we have
\begin{equation}\label{valeurs_Frob}
\begin{array}{rcl}
\Harm_{1}(\mb{w};\mb{q})\hskip-8pt&=&\hskip-8ptS_{1}(\mb{w})\\[4pt]
\Harm_{2}(\mb{w};\mb{q})\hskip-8pt&=&\hskip-8ptS_{2}(\mb{w})+s_{1}S_{11}(\mb{w})\\[4pt]
\Harm_{3}(\mb{w};\mb{q})\hskip-8pt&=&\hskip-8ptS_{3}(\mb{w})+ ( s_{1}+s_{2} )\, S_{21}(\mb{w})+ ( s_{11}+s_{3} )\, S_{111}(\mb{w})\\[4pt]
\Harm_{4}(\mb{w};\mb{q})\hskip-8pt&=&\hskip-8ptS_{4}(\mb{w})+ ( s_{1}+s_{2}+s_{3} )\, S_{31}(\mb{w})+ ( s_{2}+s_{21}+s_{4} )\, S_{22}(\mb{w})\\
         \hskip-8pt &  &\hskip-8pt+( s_{11}+s_{3}+s_{21}+s_{4}+s_{31}+s_{5} )\, S_{211}(\mb{w})\\
          \hskip-8pt &  &\hskip-8pt+ (s_{111}+ s_{31}+s_{41}+s_{6} )\, S_{1111}(\mb{w})\\[4pt]
\Harm_{5}(\mb{w};\mb{q})\hskip-8pt&=&\hskip-8ptS_{5}(\mb{w})+ ( s_{1}+s_{2}+s_{3}+s_{4} )\, S_{41}(\mb{w})\\[4pt]
          \hskip-8pt &  &\hskip-8pt+  (s_{2} +s_{3}+s_{21}+s_{4}+s_{31}+s_{22}+s_{5}+s_{41}+ s_{6})\, S_{32}(\mb{w})\\[4pt]
          \hskip-8pt &  &\hskip-8pt+ ( s_{11}+s_{3}+s_{21}+s_{4}+2\,s_{31}+2\,s_{5}+s_{41}\\
          \hskip-8pt &  &\hskip-8pt \qquad  \qquad\qquad\qquad\qquad+s_{32}+s_{6}+s_{51}+ s_{7} )\, S_{311}(\mb{w})\\[4pt]
          \hskip-8pt &  &\hskip-8pt+ ( s_{21}+s_{4}+s_{31}+s_{22}+s_{211}+s_{5}+2\,s_{41}+s_{32}\\
          \hskip-8pt &  &\hskip-8pt \qquad  \qquad\qquad+s_{311}+s_{6}+2\,s_{51}+s_{42}+s_{7}+s_{61}+s_{8} )\, S_{221}(\mb{w})\\[4pt]
          \hskip-8pt &  &\hskip-8pt+ ( s_{111}+s_{31}+s_{211}+2\,s_{41}+s_{32}+s_{311}+s_{6}+2\,s_{51}\\
         \hskip-8pt &  &\hskip-8pt \qquad  \qquad+s_{42}+s_{411}+s_{33}+s_{7}+2\,s_{61}+s_{52}+s_{8}+s_{71}+s_{9})\, S_{2111}(\mb{w})\\[4pt]
          \hskip-8pt &  &\hskip-8pt+ ( s_{311}+s_{42}+s_{411}+s_{61}+s_{511}+s_{43}\\
          \hskip-8pt &  &\hskip-8pt\qquad  \qquad\qquad\qquad\qquad+s_{71}+s_{62}+s_{81} +s_{10})\, S_{11111}(\mb{w})
\end{array}
\end{equation}
In these expressions,  the Schur functions $s_\mu=s_\mu(\mb{q})$ are to be evaluated in the three parameters $\mb{q}=q_1,q_2,q_3$. 
The Schur function $S_\lambda(\mb{w})$ act as placeholder for irreducible representations of $\S_n$.
Thus the multiplicities of the irreducible representations of $\S_n$ in $\Harm_{n}$, classified by the partitions $\mu$ of $n$,  appear as the coefficients of $S_{\mu}(\mb{w})$ through the evaluation of~\pref{valeurs_Frob}  at $q_1=1$, $q_2=1$ and $q_3=1$.
In particular, the value obtained this way, as  coefficient of $S_{11\cdots1}(\mb{w})$, agree formula~\pref{dim_An}.
Furthermore, the expressions in \pref{valeurs_Hilbert} may be directly calculated from  \pref{valeurs_Frob} by replacing each $S_\lambda(\mb{w})$
by the number $f_\mu$ of standard tableaux of shape $\mu$ (see~\cite{livre}).

\subsection*{Generalized harmonics} Using notions further discussed in Section~\ref{invariant}, all of these considerations generalize to the spaces of \defn{higher diagonal harmonics}
\begin{equation}\label{def_higher}
   \bleu{\Harm_{n}^{(r)}:=(\Altid_n^{r-1})\cap (\Altid_n^{r-1}\Ideal_n)^\perp},
 \end{equation}
 where \bleu{$\Altid_n$} (resp.  \bleu{$\Ideal_n$}) stands  for the ideal generated by \defn{alternating} polynomials (resp. \defn{invariant}) polynomials  in $\R$, and the orthogonal complement is with respect to the $\S_n$-invariant scalar product defined by
\begin{equation}\label{def_scalar}
   \bleu{\scalar{X^A}{X^B}:=\begin{cases}
     A!& \text{if}\ A=B, \\
      0 & \text{otherwise}.
\end{cases}}
\end{equation}
Here $A!$ stands for the product of factorials:
    \begin{displaymath}\bleu{A!:=a_1!\cdots a_n!\cdot b_1!\cdots b_n!\cdot c_1!\cdots c_n!}.\end{displaymath}
When $r$ is even, the diagonal action on $\Harm_{n}^{(r)}$ is twisted by the sign representation. 
 For $n=2$, the space $\Harm_{2}^{(r)}$ decomposes as a direct sum of its isotypic \defn{invariant} component, $\bleu{\Inv_2^{(r)}}$, and  \defn{alternant} component, $\bleu{\Alt_2^{(r)}}$:
 \begin{equation}
      \bleu{\Harm_{2}^{(r)}=\Inv_2^{(r)}\oplus \Alt_2^{(r)}},
 \end{equation} 
 which afford the respective bases
\begin{eqnarray*}
    \bleu{\Base\Inv_2^{(r)}}&=&\bleu{\{ (x_2-x_1)^a(y_2-y_1)^b(z_2-z_1)^c\ |\ a+b+c=r-1 \}},\qquad {\rm and}\\
    \bleu{\Base\Alt_2^{(r)}}&=&\bleu{\{ (x_2-x_1)^a(y_2-y_1)^b(z_2-z_1)^c\ |\ a+b+c=r \}}. 
  \end{eqnarray*}
It follows that
    \begin{displaymath}\bleu{\Harm_{2}^{(r)}(\mb{w};\mb{q})=s_{r-1} S_2(\mb{w}) + s_{r} S_{11}(\mb{w})},\end{displaymath}
from which we get the Hilbert series 
  \begin{displaymath}\bleu{\Harm_{2}^{(r)}(\mb{q})=s_{r-1} + s_{r} },\quad \bleu{\Inv_2^{(r)}(\mb{q})=s_{r-1}},\quad {\rm and}\quad  \bleu{\Alt_2^{(r)}(\mb{q})=s_{r}}.\end{displaymath}
In particular,
  \begin{displaymath}\bleu{\dim \Harm_{2}^{(r)} = (r+1)^2 },\quad \bleu{\dim \Inv_{2}^{(r)} = \frac{r(r+1)}{2}},\quad {\rm and}\quad  \bleu{\dim \Alt_{2}^{(r)} = \frac{(r+1)(r+2)}{2}}.\end{displaymath}
 which agree with formulas~\pref{dim_Hrn} and \pref{dim_Altr}.

The above formulas, which clearly generalize \pref{dim_Hn} and \pref{dim_An}, are both encompassed in the following generalization of a Conjecture of~\cite{HHLRU}, in conjunction with a second conjecture below.
\begin{conjecture}\label{conj1}
The graded Frobenius characteristic of   $\Harm_n^{(r)}$    is given by the formula
 \begin{equation}\label{formule_frob_r} 
   \bleu{\Harm_n^{(r)}(\mb{w};q_1,q_2,1)=\sum_{\beta\in \Dyck{r}{n}} \sum_{ f\in\Park(\beta)}  i^{(r)}_\beta(q_1)\, q_2^{\dinv(f)} Q_{\co(f)}(\mb{w})},
\end{equation}
where $\beta$ run over the set of $r$-Dyck paths, on which we consider the $r$-Tamari order \hbox{\rm (see Section~\ref{tamari})}, and $\bleu{\Park(\beta)}$ denotes the set of ``parking functions'' of ``shape'' $\beta$.  
\end{conjecture}
To finish parsing the right-hand side of \pref{formule_frob_r}, recall from~\cite{HHLRU} that $Q_{\co(f)}(\mb{w})$ stands for the fundamental quasi-symmetric polynomials indexed by  the composition  $\co(f)$ associated to a parking function $f$ as described in Section~\ref{tamari}. We also recall in Section~\ref{tamari} other relevant notions pertaining to Dyck paths, parking functions and the Tamari poset, including the definition of the ``$\dinv$'' statistic.
For a given $r$-Dyck path $\beta$, let us denote by $\bleu{ i^{(r)}_\beta}(q) $ the polynomial 
    \begin{displaymath}\bleu{ i^{(r)}_\beta(q) :=\sum_{\alpha\leq \beta} q^{d(\alpha,\beta)}},\end{displaymath} 
that $q$-enumerates the $r$-Dyck paths $\alpha$ that are smaller than $\beta$ in the $r$-Tamari order. Such paths are weighted by $q^{d(\alpha,\beta)}$, where $\bleu{d(\alpha,\beta)}$ stands for the of    \defn{length of  the longest chain} going from $\alpha$ to $\beta$ in $\Dyck{r}{n}$. 
We write $i^{(r)}_\beta$ for $i^{(r)}_\beta(1)$.
Since $d(\alpha,\beta)=0$ if and only if $\alpha=\beta$, we have $ i^{(r)}_\beta(0)=1$. 
Thus Conjecture~\ref{conj1} agrees with the shuffle conjecture of~\cite{HHLRU} at $q_1=0$. It also follows from results of~\cite{HHLRU}  that the specialization at   $q_1=q_2=1$ of formula \pref{formule_frob_r} gives
 \begin{equation}\label{formule_frob_111} 
   \bleu{\Harm_n^{(r)}(\mb{w};1,1,1)= \sum_{\beta\in\Dyck{r}{n}}     i^{(r)}_\beta\, e_{\co(\beta)}}.
 \end{equation}
 Here $e_{r_1\cdots r_k}:=e_{r_1}\cdots e_{r_k}$, where $e_r$ stands for the usual $r^{\rm th}$ elementary symmetric function (see \cite{macdonald}).
It may be worth underlying that the right hand side of this last equation is simply, up to twisting by the sign,  the Frobenius characteristic of the action of the symmetric group on the set of pairs $(f,\alpha)$, where $f$ is a $r$-parking function and $\alpha$ is a $r$-Dyck path lying below the shape of $f$. See section~\ref{tamari} for more on this.
Calculations suggest that we have the following explicit formula for this Frobenius characteristic.
\begin{conjecture}\label{conj2}
 \begin{equation}\label{formule_frob_p} 
   \bleu{\Harm_n^{(r)}(\mb{w};1,1,1)}=\sum_{\lambda \vdash n} \bleu{(-1)^{n-\ell(\lambda)}(r\, n+1)^{\ell(\lambda)-2}
                      \prod_{k\in\lambda} \binom{(r+1)\,k}{k}}\,\frac{1}{z_\lambda} p_\lambda(\mb{w}).
\end{equation}
\end{conjecture}
Recall that, for a partition $\lambda$ having $d_i$ parts of size $i$, it is usual to denote by $z_\lambda$ the integer
    \begin{displaymath}\bleu{z_\lambda:=1^{d_1}d_1!\,2^{d_2}d_2!\,\cdots\, n^{d_n}d_n!}.\end{displaymath}
It is interesting to note that the special case $r=1$ of this last formula has also been conjectured independtely by Loktev (see~\cite{loktev}).
It follows from~\pref{formule_frob_p} that we have an explicit formula for the enumeration of pairs $(f,\alpha)$, as above, that are fixed points for the action of a permutation
having given cycle type. In generating function format, this may be stated as
 \begin{equation}\label{polya}
   \bleu{\sum_{n\geq 0} \sum_{\beta\in\Dyck{r}{n}}     i^{(r)}_\beta\, h_{\co(\beta)}=
        \frac{1}{(rn+1)^2} \exp\left(\sum_{k\geq 1} (r\,n+1)\binom{(r+1)\,k}{k}\,p_k(\mb{w})/k \right)},
\end{equation}
where the $e$-basis has been turned into the $h$-basis by the removal of signs.

This also implies  that the dimension of $\Harm_n^{(r)}$ may be described in two different ways, giving rise to the following elegant  formula
\begin{equation}\label{dimension}
  \bleu{ \sum_{\beta\in\Dyck{r}{n}}     i^{(r)}_\beta\,  \binom{n}{\co(\beta)}=(r+1)^n\,(r n+1)^{n-2}}.
\end{equation}
A first draft of this paper prompted a collaborative effort between  Bous\-quet-M\'elou, Chapuy, and the second author, to give a direct proof of this new combinatorial identity  (see~\cite{chapuy}).

On another note, we will recast formula~\pref{formule_frob_p} using the following calculations. We start with establishing the formal power series identity
\begin{equation}\label{fonct_gen}
   \bleu{\exp\left(-a\,\sum_{k\geq 1} \binom{b\,k}{k}\frac{(-t)^k}{k}\right) =
  1+  \sum_{n\geq 1} \Big((a-n)\,b+1\Big)^{(n-1)}\ {a\,b}\, \frac{t^n}{n!}},
   \end{equation}
where we use the Pochhammer symbol notation
    $$(u)^{(m)}:=u\,(u+1)\,(u+2)\,\cdots \, (u+m-1).$$
To check that the above equality hold, we may proceed as follows. Let us denote by $Z=Z(t)$  the series corresponding to the right-hand side of \pref{fonct_gen}, at $a=1$. Then Equation~\pref{fonct_gen} is equivalent to
\begin{equation}\label{fonct_gen_dlog}
   \bleu{\sum_{k\geq 1} \binom{b\,k}{k}\, (-t)^{k-1}=
  \frac{Z'}{Z}},
   \end{equation}
On the other hand, for $b\in\N$, we see directly that the series $Z$  is algebraic and satisfies the  equation 
     \begin{equation}\label{eq_fonct_gen}
         \bleu{Z^{b-1}=(Z-t)^{b}},
    \end{equation}
  with initial condition $Z(0)=1$.
But \pref{eq_fonct_gen} implies that we have
     \begin{equation}\label{diffeq_fonct_gen}
         \bleu{\frac{Z'}{Z}=\frac{b}{Z+(b-1)\,t}},
    \end{equation}
and identity~\pref{fonct_gen_dlog} follows from this.

We now use Formula~\pref{fonct_gen}, setting $a=r\,n+1$ and $b=r+1$, in conjunction with the symmetric function identity
\begin{equation}
    \exp\left(-\sum_{k\geq 1} \bleu{p_k(\mb{u})} \frac{(-t)^k}{k}\right)=\sum_{n\geq 0} \bleu{e_n(\mb{u})} t^n,\
 \end{equation}
to ``evaluate''  the elementary symmetric function $e_n(\mb{u})$ for the plethystic substitution (see~\cite{livre})  
\begin{equation}\label{subst}
   \bleu{p_k(\mb{u})} \mapsto \bleu{(r\,n+1) \binom{(r+1)\,k}{k}}.
 \end{equation}
Recall that this means that we should expand $e_n(\mb{u})$ in terms of the power-sums $p_k(\mb{u})$, and then apply~\pref{subst}. In view of \pref{fonct_gen}, 
we get that
\begin{equation}\label{subste}
   \bleu{e_n(\mb{u})} \mapsto \bleu{\frac{r\,n+1}{r n-k+1}
      \binom{r((r+1)n-k+1)}{k}}.
 \end{equation}
This calculation allows us to reformulate \pref{formule_frob_111} in terms of the monomial symmetric functions $m_\lambda(\mb{w})$ by applying the substitution~\pref{subst} to the Cauchy identity (see \cite[page 65]{macdonald}):
 \begin{equation}\label{cauchy}
    \sum_{\lambda\vdash n} (-1)^{n-\ell(\lambda)}\bleu{p_\lambda(\mb{u})}\,\frac{1}{z_\lambda} p_\lambda(\mb{w})=\sum_{\lambda\vdash n}  \bleu{e_\lambda(\mb{u})}\, m_\lambda(\mb{w}).
  \end{equation}
We thus get
  \begin{equation}
  \begin{array}{l}
\displaystyle \sum_{\lambda \vdash n} (-1)^{n-\ell(\lambda)}\bleu{(r\, n+1)^{\ell(\lambda)}
                      \prod_{k\in\lambda} \binom{(r+1)\,k}{k}}\,\frac{1}{z_\lambda} p_\lambda(\mb{w}) =\\[20pt]
 \hskip2cm\displaystyle\sum_{\lambda\vdash n}  \bleu{(r n+1)^{\ell(\lambda)} \prod_{k\in \lambda} \frac{1}{r n-k+1}
      \binom{r((r+1)n-k+1)}{k}}\, m_\lambda(\mb{w}),
  \end{array}
  \end{equation}
simply by applying the substitution~\pref{subst} in the left-hand side of \pref{cauchy}, and \pref{subste} in the right-hand side.
  
Formula~\pref{formule_frob_p} can thus be rewritten in the form
 \begin{equation}\label{formule_frob_m}
  \bleu{{\Harm_n^{(r)}(\mb{w};1,1,1)}=\sum_{\lambda\vdash n}  (r n+1)^{\ell(\lambda)-2} \prod_{k\in \lambda} \frac{1}{r n-k+1}
      \binom{r((r+1)n-k+1)}{k}}\, m_\lambda(\mb{w}).
\end{equation}
In particular, the multiplicity of the alternating representation in $\Harm_n^{(r)}$ is the coefficient of $m_{11\cdots 1}(\mb{w})$ in this last expression. In view of \pref{formule_frob_111}, this suggested that we have the following generalization of a formula of Chapoton (see \cite{chapoton}, in the case $r=1$) for the number of intervals in the $r$-Tamari poset: 
 \begin{equation}\label{intervalles_m_Tamari}
  \bleu{\sum_{\beta\in\Dyck{r}{n}}    i^{(r)}_\beta = {\frac {( r+1)}{ n\,( r n+1) }} \binom{( r+1) ^{2}\,n+r}{ n-1} }.
\end{equation}
A few months before this writing, the first author conjectured that this last equality should hold. This new combinatorial identity has since been shown to be true in~\cite{melou}. 

In view of Proposition~\ref{eqalite_inv_alt},  the above formulas and conjectures  imply that the following combinatorial identity should also hold
    \begin{equation}\label{partie_triviale}
  \bleu{ \sum_{\beta\in\Dyck{r}{n}} i^{(r)}_\beta\  \charac(\co(\beta)=11\cdots 1)  \ =
       \sum_{\beta\in\Dyck{r-1}{n}} i^{(r-1)}_\beta  },
\end{equation}
where $\charac(P)$ is equal to $1$ if $P$ is true, and $0$ otherwise.
A direct bijective combinatorial proof of this identity seems to follow from the approach in~\cite{melou}. Moreover, it appears that \pref{polya}
may be amenable to a Polya theory extension of the approach of~\cite{chapuy} and~\cite{melou}. In particular, this will reduce Conjecture~\ref{conj2} to Conjecture~\ref{conj1},
as well as giving a joint proof of \pref{dimension} and \pref{intervalles_m_Tamari}.

\section{Trivariate diagonal invariants and alternants}\label{invariant}
\subsection*{Invariants and scalar product} For our purpose, some interesting properties of trivariate diagonal invariants and alternants need to be recalled. 
Much of what is discussed here is a straightforward generalization to the context of three sets of variables of results of \cite{garsia_haiman} for the bivariate case. They are included here to help make our presentation understandable on its own. Omitted proofs are entirely similar to the bivariate case.
A particular case of a result of Weyl, in~\cite{weyl}, establishes that the subring of trivariate diagonal invariants is generated by the \defn{polarized power sums}
  \begin{equation}\label{def_power_sum}
      \bleu{P_\alpha(X):=\sum_{j=1}^n X_j^\alpha},
  \end{equation}
for $\alpha=(a,b,c)$ running thought the set of all $\N$-vectors such that $1\leq |\alpha|\leq n$. 

To make better sense out of Definition~\pref{def_higher}, we need to discuss a few notions concerning trivariate diagonal alternants.
Let us first observe that definition~\pref{def_scalar} is equivalent to 
\begin{equation}\label{def_scalar2}
    \bleu{ \scalar{f(X)}{g(X)}=f(\partial X)g(X)\big|_{X=0}},
 \end{equation}
where $\bleu{p(X)\big|_{X=0}:=p(0)}$ is the \defn{constant term} of $p(X)$. 

It is a well known fact that, with respect to this scalar product, the orthogonal complement $I^\perp$, of an ideal  $I=(f_1,\ldots, f_N)$, is simply  the set of polynomial zeros of the differential operators $f_i(\partial X)$. Indeed, let $g(X)$ be in $I^\perp$, and consider the leading term (for any suitable term order) of
    \begin{displaymath} f_i(\partial X)g(X) = c\, X^A+ \ldots, \end{displaymath}
 if any.
Since $X^A\,f_i(X)$ lies in $I$, we must have $\scalar{ X^A\,f_i(X)}{g(X)} =0$, but this means precisely that $c=0$. Hence we must have $f_i(\partial X)G(X)=0$. In particular, Weyl's result implies that 
      \begin{displaymath}\Ideal_n^\perp =\Harm_n^{(1)}=\Harm_n.\end{displaymath} 
To show that the spaces we consider are finite dimensional, the following is very useful.      
\begin{proposition}\label{eqalite_inv_alt}
      All monomials $X^A$, of total degree larger than $\binom{n}{2}$ lie in $\Ideal_n$.
\end{proposition}

      
\subsection*{Alternants} A set of linear generators for trivariate alternants is easily obtained by applying the \defn{antisymetrization operator}
   \begin{equation}\label{defn_reynold}
        \bleu{R^\pm(X^A) :=\sum_{\sigma\in\S_n} \sign(\sigma) \, \sigma\cdot X^A},
     \end{equation}   
 to any monomial $X^A$, with $A$ a $3\times n$ matrix of non-negative integers.
Evidently the product of an alternant by an invariant is yet again an alternant.
We get a linear basis $\Delta_A(X):=\{R^\pm(X^A)\}_A$, of the module of alternants over the ring of invariants, by choosing matrices $A$ having all their columns different and appearing in decreasing lexicographic order (from left to right). Any symmetric operator in the ring $\C[X,\partial X]$ sends alternants to alternants. In particular, this is the case for the operators $E_{\mb{u}\mb{v}}^{(k)}$, and $P_\alpha(\partial X)$. A degree argument shows that  $P_\alpha(\partial X)\Delta_A(X)=0$, whenever $A$ is an order-ideal under \defn{component} comparison of columns:
    \begin{displaymath}\bleu{B\leq C,\qquad {\rm iff}\qquad \forall i\ (b_i\leq c_i)}.\end{displaymath}
 For instance, this is the case for the Vandermonde determinants. 
  Moreover, a simple direct calculation of the commutator of the operators  $E_{\mb{y}\mb{x}}^{(k)}$ and $P_\alpha(\partial X)$, with $\alpha=(a,b,c)$,  gives that
 \begin{equation}\label{commutator}
     \bleu{[P_\alpha(\partial X),E_{\mb{y}\mb{x}}^{(k)}]=a\, P_\beta(\partial X)},
 \end{equation} 
where \bleu{$\beta=\alpha+(k,-1,0)$}.
 It is evident that similar results holds for other possible choices of variable sets in $E_{\mb{u}\mb{v}}^{(k)}$. 
Putting all this together, we deduce that
 \begin{proposition}\label{operator_calculation}
    If $A\in \N^{3\times n}$ is an order-ideal, then 
     \begin{equation}\label{operator}
        \bleu{(\partial X^A)\Omega ( \Delta_A(X))\in \Harm_n},
     \end{equation}
for any $A$, and any composition $\Omega$ of the operators $E_{\mb{u}\mb{v}}^{(k)}$.
  \end{proposition}
The ``Operator Theorem'' of~\cite{haimanvanishing}  states that Proposition~\ref{operator_calculation} entirely characterizes the space $\Harm_n$ in the bivariate case. More precisely, we have the following.
\begin{thm}[Haiman]
   The space of diagonal harmonics for the bivariate case is the smallest vector space containing the Vandermonde determinant $\Delta_n(\mb{x})$, which is closed under taking partial derivatives and applications of  the operators $E_{\mb{y}\mb{x}}^{(k)}$. 
 \end{thm}
 For the trivariate case, experiments suggest that the analogous statement should hold. However,  the methods employed by Haiman to settle the bivariate case do not seem to generalize to the trivariate case. 

Another result along these lines is that, for the bivariate case, we can calculate the entire space $\Alt_n^{(r)}$ (resp. $\Inv_n^{(r)}$) by successive applications of the operators  $E_{\mb{y}\mb{x}}^{(k)}$, starting with the $r^{\rm th}$-power (resp. $(r-1)^{\rm th}$-power) of the Vandermonde determinant $\Delta_n(\mb{x})$. However, we do not know how to explicitly construct a basis in this manner. All, this also appears to hold in the trivariate case. In particular, we have that

\begin{proposition}\label{eqalite_inv_alt}
For all $n$, and all $r\geq 1$, the spaces $\Inv_n^{(r)}$ and $\Alt_n^{(r-1)}$ coincide.
\end{proposition}

\section{Action of the general linear group}\label{general}
To make apparent a special feature of both the Hilbert series and the (graded) Frobenius characteristic of the spaces $\Harm_n^{(r)}$, we need to consider it as a polynomial $GL_{3}$-sub-representation of $\R_n=\C[X]$, for the action
    \begin{equation}\label{def_action_GL}
         \bleu{f(X) \mapsto f(M\,X)},\qquad {\rm for}\qquad \bleu{M\in GL_{3}}.
      \end{equation}
Clearly this action commutes with the action of $S_n$. Recall that the \defn{character} of a $GL_3$ representation is the symmetric function of the parameters \begin{math} \mb{q}=q_1,q_2,q_3\end{math} obtained by calculating the trace of the linear transformation
   \begin{displaymath} \bleu{Q^*(f(X)) :=f(Q\,X)},\end{displaymath}
 where $Q$ is the diagonal matrix $Q=[q_1,q_2,q_3]$.
 Observe that a polynomial $f(X)$ is \defn{homogeneous}  of degree $d$ if and only if 
    \begin{displaymath} \bleu{f(Q\,X)=\mb{q}^d\, f(X) }.\end{displaymath}    
Since an homogenous subspace $\mathcal{V}$, of  $\R_n$, is simply a subspace that affords a basis $\mathcal{B}$ of homogeneous polynomials,
 it is $G_3$-invariant and its character $\mathcal{V}(\mb{q})$ is 
    \begin{displaymath} \bleu{\mathcal{V}(\mb{q})=\sum_{f\in \mathcal{B}} \mb{q}^{\deg(f)}} .
    \end{displaymath} 
 This is precisely what we have called previously the Hilbert series of $\mathcal{V}$. It is well known that the characters, of irreducible polynomial $GL_3$-representations, are the Schur polynomials (evaluated at $q_1$, $q_2$, $q_3$). Hence, $\mathcal{V}(\mb{q})$ expands as a positive integer coefficient linear combination of the $s_\mu(\mb{q})$, since the coefficients correspond to multiplicities of irreducibles. 
 
 Since the two actions commute, we get a double decomposition into irreducibles, for both $\S_n$ and $GL_3$, of any subspace that is stable under the two actions. This is the reason why the graded Frobenius characteristic of such a space expands as a positive sum of products  $s_\mu(\mb{q}) S_\lambda(\mb{w})$.

\section{The whole story for the case $n=3$}\label{nabla}
Using observations and results of~\cite{bergeron_several}, we can calculate explicitly $\Harm_3^{(r)}(\mb{w};\mb{q})$ in full generality. The crucial observation is that the only Schur function $s_\mu(\mb{q})$ that can occur in the expansion of $\Harm_3^{(r)}(\mb{w};\mb{q})$ are those indexed by partitions having at most two parts. We can thus deduce the trivariate case\footnote{In fact, this particular calculation holds for any number of set of variables.} expansion from the known bivariate case.  

Recall that, in the bivariate case we have 
$ {\Harm_3^{(r)}(\mb{w};q,t)=\nabla^r(S_{111}(\mb{w}))}$,
   where $\nabla$ is an operator characterized as follows.
As usual (see \cite{macdonald}), we denote by $\lambda\preceq\mu$ the \defn{dominance order} on partitions, and $\mu'$ stands for the \defn{conjugate} partition of $\mu$. 
The \defn{integral form Macdonald polynomials} $H_\mu(\mb{w};q,t)$, with $\mu$ a partition of $n$  form a linear basis of the ring $\Lambda$, of symmetric functions in the variables $\mb{w}$, and   are characterized by the equations 
 \begin{equation}\label{eq_Hmu}
 \begin{array}{llll}
    \mathrm{(i)}\ \bleu{\displaystyle  \langle S_\lambda(\mb{w}), H_\mu[(1-q)\,\mb{w};q,t]\rangle=0},\qquad{\rm if}\qquad  \lambda\not\succeq\mu,\\[6pt]
    \mathrm{(ii)}\ \bleu{\displaystyle   \langle S_\lambda(\mb{w}), H_\mu[(1-t)\,\mb{w};q,t]\rangle=0},\qquad{\rm if}\qquad   \lambda\not\succeq\mu',\ \mathrm{and}\\ [6pt]
    \mathrm{(iii)}\ \bleu{\displaystyle   \langle S_n(\mb{w}), H_\mu(\mb{w};q,t)\rangle=1},
      \end{array}
 \end{equation} 
 involving the usual ``Hall'' scalar product on symmetric functions (for which the Schur functions are orthonormal), and some plethystic notation (see~\cite{livre}). The linear operator $\nabla$  is characterized by the fact that it affords the Macdonald polynomials as eigenfunctions:
     \begin{displaymath}  
         \bleu{\nabla(H_\mu(\mb{w};q,t))=H_\mu(\mb{w};q,t)\, \prod_{(a,b)\in \mu} q^at^b}.
      \end{displaymath}  
This presentation may seem a bit terse, but for our purpose we only need the case $n=3$ for which an explicit description of $\nabla$ is given below. Using this description, we will infer the graded Frobenius of $\Harm_3^{(r)}$ for the trivariate case, out of the simple knowledge of the bivariate case. To this end, we  consider the linear operator  corresponding to the matrix
\begin{equation}\label{formule_nabla}
\nabla=\begin{pmatrix}
   0&s_{{22}}&s_{{32}}\\ \noalign{\medskip}0
&-s_{{21}}&-s_{{31}}\\ \noalign{\medskip}1&s_{{1}}+s_{{2}}&s_{{11}}
+s_{{3}}
\end{pmatrix},
\end{equation}
with respect to the basis  $\{S_3(\mb{w}),S_{21}(\mb{w}),S_{111}(\mb{w})\}$. Then\footnote{All calculations are to be done over the ring $\Lambda^{(2)}(\mb{q})$ of symmetric functions modulo the ideal generated by the $s_\mu$ for which  $\mu$ has at least three parts.},
\begin{proposition}\label{prop_itnabla}
For the trivariate case, we have the explicit formula
  \begin{equation}\label{eq_itnabla}
      \bleu{\Harm_3^{(r)}(\mb{w};\mb{q})=\nabla^r(S_{111}(\mb{w}))}.
  \end{equation}
 \end{proposition}
We may calculate from \pref{formule_nabla} that
\begin{lemma}\label{prop_charpoly}
  The characteristic polynomial of the operator $\nabla$, over the ring $\Lambda^{(2)}(\mb{q})$, is
  \begin{equation}\label{form_charpoly}
    \bleu{  \nabla^3- ( s_3-s_{21}+s_{11} ) \nabla^{2}+ ( s_{41}+s_{33}-s_{32} ) \nabla-s_{44}}.
   \end{equation}
\end{lemma}
One may observe that $s_{11}(\mb{q})=q_1q_2+q_1q_3+q_2q_3$ is one of the eigenvalues of $\nabla$. The associated eigenspace is spanned by the
symmetric function
    $$S_3(\mb{w}) + (q_1+q_2+q_3)\,S_{21}(\mb{w})+(q_1q_2+q_1q_3+q_2q_3)\,S_{111}(\mb{w}),$$
which is a three parameter version of a Macdonald polynomial. Wether one can develop a  theory of these polynomials for all $n$ remains an open question. One would expect these to appear as eigenfunctions for a suitable generalization of $\nabla$.

In view Proposition~\ref{prop_itnabla}, Lemma~\ref{prop_charpoly}  leads to a recursive approach for the calculation of the graded Frobenius characteristic of any of the spaces $\Harm_3^{(r)}$, based on the initial values
\begin{eqnarray*}
     \bleu{\Harm_3^{(0)}(\mb{w};\mb{q})}&=& \bleu{S_{111}(\mb{w})},\\
     \bleu{\Harm_3^{(1)}(\mb{w};\mb{q})}&=& \bleu{S_3(\mb{w})+(s_{2} +s_{1} )\, S_{21}(\mb{w})+(s_{3} +s_{11})\, S_{111}(\mb{w})},\\
     \bleu{\Harm_3^{(2)}(\mb{w};\mb{q})}&=& \bleu{( s_{11}+s_{3} ) S_{3}(\mb{w})+ ( s_{21}+s_{4}+s_{31}+s_{5} ) S_{21}(\mb{w})}\\
              &&\qquad \bleu{+   ( s_{2,2}+s_{41}+s_{6}) S_{111}(\mb{w})}.
 \end{eqnarray*}
 Calculations are to be done in the ring $\Lambda^{(2)}(\mb{q})$ of symmetric functions spanned by Schur functions indexed by partitions that have at most two parts. In order to describe the resulting formulas, let us introduce the notation
     \begin{equation}
        \bleu{T_{nk}:= s_n+s_{n-2,1}+\,\ldots\,+ s_{n-2\,k,k}}.
     \end{equation}
After calculations, we get the following.
\begin{proposition}
The space $\Harm_3^{(r)}$ affords the explicit graded Frobenius character formula
\begin{equation}\label{Frob_explicit3}
      \bleu{\Harm_3^{(r)}(\mb{w};\mb{q})= T_{3\,(r-1),r-1}\,S_3(\mb{w})+ (T_{3\,r-2,r-1}+T_{3\,r-1,r-1 })\,S_{21}(\mb{w})+ T_{3\,r,r}\,S_{111} (\mb{w})   }.
\end{equation}
In particular, the Hilbert series of the alternating part $\Alt_3^{(r)}$ of this space appears as the coefficient of $S_{111}(\mb{w})$ in this formula, i.e.:
 \begin{equation}
      \bleu{\Alt_3^{(r)}(\mb{q})= s_{3\,r} +s_{3\,r-2,1} +\,\ldots\, + s_{r,r}  }.
\end{equation}   
\end{proposition}
For any $a\geq b\geq 0$, the Schur function $s_{ab}=s_{ab}(q_1,q_2,q_3)$ specializes to
    $$\bleu{\frac{1}{2} (b+1)(a+2) (a-b+1)},$$
when one sets $q_1=q_2=q_3=1$.
From this, it follows that
    $$\bleu{T_{nk}(1,1,1)=\frac{1}{12}\, \left( k+1 \right)  \left( k+2 \right)  \left( {k}^{2}- \left( 
13+10\,n \right) k+3\, \left( n+1 \right)  \left( n+2 \right) 
 \right) 
}$$
We may then check directly that formula \pref{Frob_explicit3} specializes at $q_1=q_2=q_3=1$ to the case $n=3$ of  \pref{formule_frob_p}.

\section{Dyck paths, Tamari posets, and parking functions}\label{tamari}
\subsection*{Generalized Dyck paths and Tamari posets} 
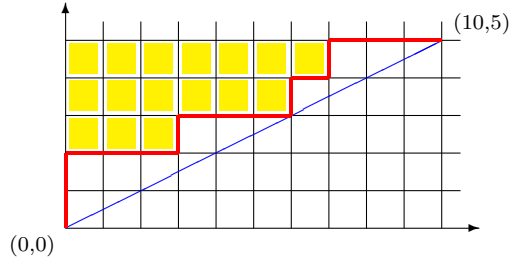
\begin{figure}[ht]
\setlength{\unitlength}{5mm}
\setlength{\carrelength}{4mm}
\def\jcarre{\put(.1,-.48){\jaune{\linethickness{\carrelength}\line(1,0){.75}}}}
\begin{center}
\begin{picture}(11,7)(0,0)
\multiput(0,3)(1,0){3}{\jcarre}
\multiput(0,4)(1,0){6}{\jcarre}
\multiput(0,5)(1,0){7}{\jcarre}
\multiput(0,1)(0,1){5}{\line(1,0){10.5}}
\multiput(1,0)(1,0){10}{\line(0,1){5.5}}
 \put(0,0){\vector(1,0){11}}
 \put(0,0){\vector(0,1){6}}
\put(-1.5,-.6){$\scriptstyle(0,0)$}
 \put(0,0){\bleu{\line(2,1){10}}}
  \linethickness{.5mm}
\put(0,0){\rouge{\line(0,1){1}}}
\put(0,1){\rouge{\line(0,1){1}}}
\put(0,2){\rouge{\line(1,0){3}}}
\put(3,2){\rouge{\line(0,1){1}}}
\put(3,3){\rouge{\line(1,0){3}}}
\put(6,3){\rouge{\line(0,1){1}}}
\put(6,4){\rouge{\line(1,0){1}}}
\put(7,4){\rouge{\line(0,1){1}}}
\put(7,5){\rouge{\line(1,0){3}}}
\put(10.3,5.3){$\scriptstyle(10,5)$}
\end{picture}\end{center}
\caption{The $2$-Dyck path encodes as $00367$.}
\label{fig1}
\end{figure}

Recall that a $r$-Dyck path of height $n$ is a path in $\N\times \N$, that starts at $(0,0)$, ends at $(rn,n)$, and stays above the line $x=ry$. For instance, we have the $2$-Dyck path of Figure~\ref{fig1}.
Such a path may be entirely described by a sequence $\alpha=(a_1,a_2,\ldots,a_n)$, where $(a_i,i)$ corresponds to the leftmost point lying on the path at height $i$. 
I other words, $a_i$ is the number of ``boxes'' lying to the left of the path at height $i$.
In this manner, $\alpha$ is a \defn{$r$-Dyck path} if and only if
\begin{enumerate}\itemsep=4pt
   \item[(1)] $\bleu{0\leq a_1\leq a_2\leq \ldots \leq a_n}$,
   \item[(2)] for each $i$, we have $\bleu{a_i\leq r(i-1)}$.
\end{enumerate}
We denote by $\bleu{\Dyck{r}{n}}$ the  set of $r$-Dyck paths of height $n$. For example,
  \begin{displaymath} \Dyck{2}{3}=\left\{ 000,001,002,003,004,011,012,013,014,022,023,024 \right\} \end{displaymath}
It is well known (see \cite{duchon,fuss})  that the number of $r$-Dyck paths of height $n$ is given by the Fuss-Catalan number:
   \begin{displaymath}\bleu{\# \Dyck{r}{n}=\frac{1}{rn+1}\binom{(r+1)\,n}{n}}.\end{displaymath}

We say that $(a_i,a_{i+1},\ldots, a_k)$ is a \defn{primitive subsequence}, of a $r$-Dyck path $(a_1,\ldots,a_n)$,  if  
\begin{enumerate}\itemsep=4pt
   \item[(1)] $\bleu{a_j-a_i< r(j-i)}$ for each $i< j\leq k$, and 
   \item[(2)] either $\bleu{k=n}$, or $\bleu{a_{k+1}-a_i\geq r(k+1-i)}$. 
\end{enumerate}
   For each $i$, there is a unique such primitive subsequence. It corresponds to the portion of the $r$-Dyck path that starts at $(a_i,i)$ and ends at the ``first return'' of the path to the line of slope $1/r$ passing through the point $(a_i,i)$.
Whenever $i$  is such that
$a_{i-1}<a_i$, we set $\bleu{\alpha\leq \beta}$, where 
    \begin{displaymath}\bleu{\beta:=(a_1,\ldots,a_{i-1},a_i-1,\ldots,a_k-1,a_{k+1},\ldots,a_n)},\end{displaymath}
with $(a_i,\ldots,a_k)$ equal to the primitive subsequence starting at $a_i$. The  \defn{$r$-Tamari poset}\footnote{As far as we know this poset has not been considered before.}, on the underlying set $\Dyck{r}{n}$, is the order obtained as the reflexive and transitive closure of this covering relation $\alpha\leq \beta$. Its largest element is the path $00\cdots 0$, and its smallest is the path for which $a_i=r\,(i-1)$.

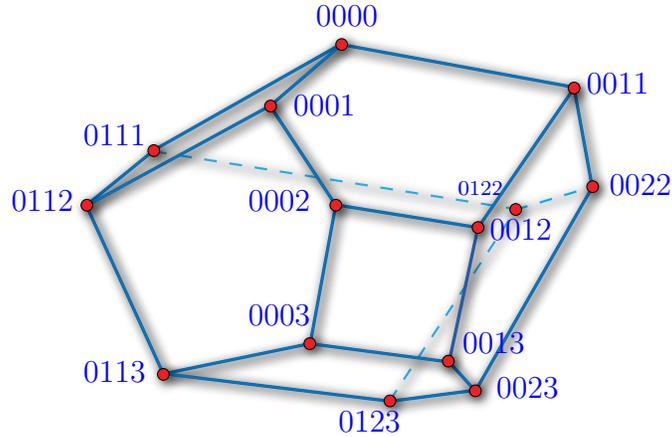
\begin{figure}[ht]
   \setlength{\unitlength}{3mm}
\begin{center}
\begin{picture}(0,0)(0,0)
                \put(19,11){$\bleu{{\scriptstyle 0122}}$}  
 \end{picture}
\scalebox{.4}{\includegraphics{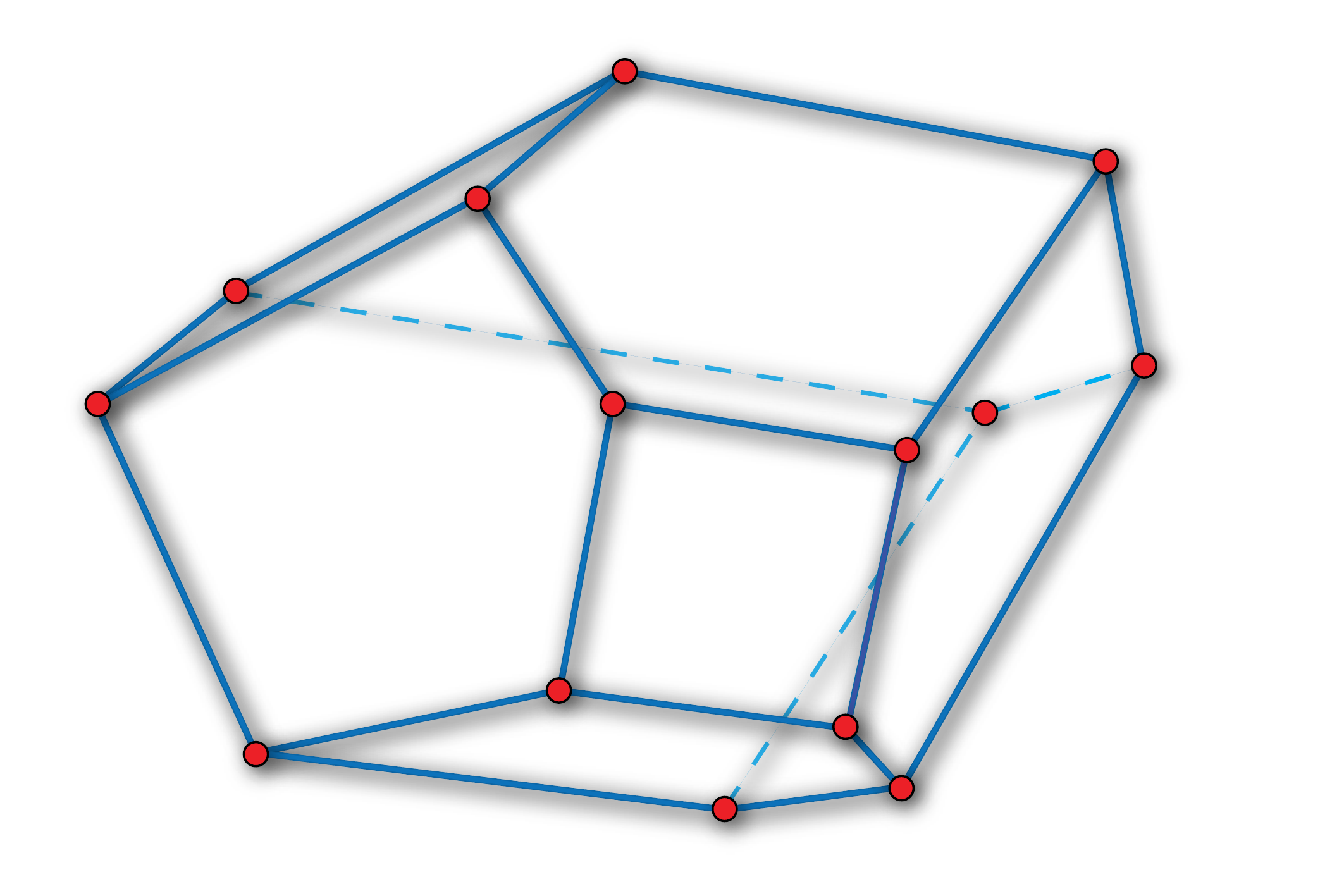}}
\begin{picture}(0,0)(0,0)
               \put(-17,18.5){$\bleu{0000}$}
       \put(-27.3,13.2){$\bleu{0111}$}        \put(-18,14.5){$\bleu{0001}$}      \put(-5,15.5){$\bleu{0011}$}
  \put(-30.5,10.3){$\bleu{0112}$}   \put(-20,10.3){$\bleu{0002}$}  \put(-9.3,9){$\bleu{0012}$}\put(-4,11){$\bleu{0022}$}
                 \put(-20,5.2){$\bleu{0003}$} 
  \put(-27.3,2.7){$\bleu{0113}$}      \put(-10.5,4){$\bleu{0013}$} \put(-9,2){$\bleu{0023}$}
              \put(-16,0.5){$\bleu{0123}$}
\end{picture}
\end{center}
\vskip-8pt
\caption{The $1$-Tamari poset for $n=4$.}
\label{tamari3}
\end{figure}

For a given $r$-Dyck path $\beta$, we have already considered the number $\bleu{ i^{(r)}_\beta} $ of $r$-paths that smaller than $\beta$ in $\Dyck{r}{n}$:
    \begin{displaymath}\bleu{ i^{(r)}_\beta :=\# \{\alpha\in \Dyck{r}{n}\ |\ \alpha\leq \beta\}},\end{displaymath} 
and  denoted by $d(\alpha,\beta)$ the length of the longest chain going from $\alpha$ to $\beta$ in $\Dyck{r}{n}$. Observe that, when $\alpha$ is the smallest element of $\Dyck{r}{n}$, $d(\alpha,\beta)$ coincides with  the usual ``area'' statistic, 
    \begin{displaymath}\bleu{{\rm area}(\beta):=r\binom{n}{2}-\sum_{i=1}^n a_i},\end{displaymath} 
on $r$-Dyck path $\beta$.
For $n=4$ and $r=1$, the Tamari poset is drawn up in Figure~\ref{tamari3}.
To a $r$-Dyck path $\beta$ we associate the composition $\co(\beta)$ whose parts are the length of runs of $1$ in $\beta$. Thus $\co(00112)=221$.  

\subsection*{Parking functions} For  a sequence of positive integers $f=(f_1,f_2,\ldots, f_n)$  (we also write $f=f_1f_2\cdots f_n$),
let 
\begin{equation}\label{alpha_beta}
   \bleu{\begin{pmatrix}
            \alpha(f)\\ \beta(f)\end{pmatrix} = 
        \begin{pmatrix}
         a_1 & a_2 & \ldots & a_n\\
         b_1 & b_2 & \ldots & b_n
         \end{pmatrix}}
    \end{equation}
  be  the \defn{lex-increasing rearrangement} of 
      \begin{displaymath}\begin{pmatrix}
         1 & 2 & \ldots & n\\
         f_1 & f_2 & \ldots & f_n
         \end{pmatrix},\end{displaymath}
  ordering first with respect to the bottom row.
This is to say that $\bleu{b_i\leq b_{i+1}}$, for all $1\leq i<n$, and that $\bleu{a_i<a_{i+1}}$ when $\bleu{b_i=b_{i+1}}$. 
One says that $f$ is a \defn{$r$-parking function}  if
  $\beta(f)$ is such that $b_k\leq  r(k-1)$, for all $k$. 
Perforce, $\beta(f)$ is thus a $r$-Dyck path. It is said to be the \defn{shape} of $f$.
  
  For example, with $r=2$, we have the following $49$  $r$-parking functions of length $3$:
\begin{displaymath}\begin {array}{ccccccc} 
000&001&002&003&004&010&011\\ 
012&013&014&020&021&022&023\\ 
024&030&031&032&040&041&042\\  
100&101&102&103&104&110&120\\  
130&140&200&201&202&203&204\\  
210&220&230&240&300&301&302\\  
310&320&400&401&402&410&420.\end {array}
\end{displaymath}
The set of $r$-parking functions of length $n$,  denoted by \bleu{$\Parking{n}$},   has cardinality $(rn+1)^{n-1}$ (see~\cite{yan}).
For a $r$-Dyck path $\beta$, we further denote by $\bleu{\Park(\beta)}$  the set of all parking functions of shape $\beta$. Observe that these are automatically $r$-parking functions.
Clearly we get a permutation action of symmetric group $\S_n$ on each of the $\Park(\beta)$ (by permutation of values). The corresponding Frobenius characteristic is easily seen to be equal to $h_{\co(\beta)}(\mb{w})$, where $\bleu{h_{c_1\cdots c_k}(\mb{w}):=h_{c_1}(\mb{w})\cdots h_{c_k}(\mb{w})}$ is the complete homogeneous symmetric function associated to the composition $\co(\beta)=c_1\cdots c_k$.
In particular, the number of parking functions of shape $\beta$ is given by the multinomial coefficient
 \begin{equation}\label{dim_parking}
     \#\Park(\beta)= \binom{n}{\co(\beta)}=\frac{n!}{c_1!c_2!\cdots c_k!}.
\end{equation}
Formula~\pref{dim_parking} reflects the fact that any representation of $\S_n$, having $h_c(\mb{w})$ as its Frobenius characteristic, must be of dimension $ \binom{n}{c}$.
From this we get the  identity
\begin{equation}\label{dim_deux}
   \bleu{(rn+1)^{n-1} = \sum_{\beta\in \Dyck{r}{n} }  \binom{n}{\co(\beta)}}
\end{equation}
Recall that it has been shown in~\cite{haimanvanishing} that the (ungraded) Frobenius characteristic of the $\mb{z}$-free component of $\Harm_n$ is equal to the character of the $\S_n$-module $\Parking{n}$, twisted by the sign representation. From the point of view of Frobenius characteristics, this twist turns $h_c(\mb{w})$ into $e_c(\mb{w})$. Hence, one consequence of the results in~\cite{haimanvanishing} may be stated as
 \begin{equation}\label{frob_parking}
     \bleu{\Harm_n^{(r)}(\mb{w};1,1,0)= \sum_{\beta\in \Dyck{r}{n}} e_{\co(\beta)}(\mb{w})}.
\end{equation}
The composition $\bleu{\co(f)}=c_1c_2\cdots c_k$, associated to a $r$-parking function $f$, is simply the composition that encodes the descents of the permutation $\alpha(f)$. 
 
 For the purpose of stating Conjecture~\ref{conj1}, and discussing potential extensions (see next section), let us recall (from \cite{HHLRU}) how to calculate the ``$\dinv$'' statistic of a parking function.  Given a $r$-parking function $f$, let $a_i$ and $b_i$ be as in \pref{alpha_beta}, and let $c_i:=r\,i-b_i$.
For $1\leq i< j\leq n$ and $0\leq d\leq r-1$, we say that a triple $(i,j,d)$ is a \bleu{$D$-inversion} of $f$ if either
\begin{itemize}
\item $c_i-c_j+d=0$ and $a_i<a_j$, or
\item $1\leq c_i-c_j+d\leq r-1$, or
\item $c_i-c_j+d=r$ and $a_i>a_j$.
\end{itemize} 
Then $\bleu{\dinv(f)}$ is simply the number of $D$-inversions of $f$. Notice that this actually depends on the value of $r$.

\section{Final remarks}
A missing component in~\pref{formule_frob_r} is a ``statistic'' $\rouge{\nu(f,\alpha)}$ that would account for the behavior of the third parameter $q_3$, with respect to pairs $(f,\alpha)$, for $f\in \Park(\beta)$ and $\alpha\leq \beta$ in $\Dyck{r}{n}$. Such a statistic would give a complete combinatorial description of $\Harm_n^{(r)}(\mb{w};q_1,q_2,q_3)$ in the form
 \begin{equation}\label{formule_frobtout} 
   \bleu{\Harm_n^{(r)}(\mb{w};q_1,q_2,q_3)= \sum_{\beta\in\Dyck{r}{n}}\sum_{ f\in\Park(\beta)} \sum_{\rouge{\alpha}\leq \beta} q_1^{d(\alpha,\beta)} q_2^{\dinv(f)} \rouge{q_3^{\nu(f,\alpha)}} Q_{\co(f)}(\mb{w})}.
\end{equation}
We would expect $\nu$ to be such that  $\nu(f,\alpha)=0$ if and only if $\alpha$ is the smallest element of $\Dyck{r}{n}$. Then, at $q_3=0$, formula \pref{formule_frobtout} would specialize precisely  to the Conjecture of Haiman  {\it et al.} in~\cite{HHLRU}. However, such a statistic still has to be found, even in the case $r=1$.

Much of what is discussed here seems to hold for generalized permutation groups $G(r,n)$. This will be explored in upcoming work.

 \subsection*{Thanks}
We thank Mark Haiman for insightful suggestions to the effect that we should extend our approach to the case of the $r$-spaces $\Harm^{(r)}_n$.


\end{document}